\documentclass[11pt]{scrartcl}

\usepackage[latin1]{inputenc}
\usepackage[T1]{fontenc}

\usepackage{lineno,hyperref}
\usepackage{amsmath, amssymb}

\newtheorem{thm}{Theorem}[section]
 
 \newtheorem{lem}[thm]{Lemma}

 \newtheorem{defn}[thm]{Definition}
 \newtheorem{rem}[thm]{Remark}
 \numberwithin{equation}{section}

\begin{document}


\title{Fractional Riesz-Hilbert transforms and fractional monogenic signals}

\author{Swanhild Bernstein\\ TU Bergakademie Freiberg \\ Institute of Applied Analysis\\ swanhild.bernstein@math.tu-freiberg.de}



\maketitle
\begin{abstract}
The fractional Hilbert transforms plays an important role in optics and signal processing. In particular the analytic signal proposed by Gabor has as a key component the Hilbert transform. The higher dimensional Hilbert transform is the Riesz-Hilbert transform which was used by Felsberg and Sommer to construct the monogenic signal. We will construct fractional and quaternionic fractional Riesz-Hilbert transforms based on a eigenvalue decomposition. We will prove properties of these transformations such as shift and scale invariance, orthogonality and the semigroup property. Based on the fractional/quaternionic fractional Riesz-Hilbert transform we construct (quaternionic) fractional monogenic signals. These signals are rotated and modulated monogenic signals.
\end{abstract}




\section{Introduction}
The Hilbert transform is well known in mathematics and in applications such as optics, signal processing and image processing. \\[1ex]
In mathematics the Hilbert transform was used to solve the Riemann-Hilbert problem which describes a relation between boundary values of analytic functions.
The first definition of the Hilbert transform is given in terms of Fourier series expansion of a function $f$ defined on the unit circle $\mathbb{T}=\{ z\in \mathbb{C}:\,|z| =1\}.$  Let $f\in L^1(\mathbb{T})$ with Fourier expansion $f(t) = \sum_{n=-\infty}^{\infty} c_n\,e^{2\pi\,int}$ then 
$$(Hf)(t):= -i\sum_{n=-\infty}^{\infty} \mathrm{sgn\,}(n)\,c_n\,e^{2\pi\,int}. $$
Specifically, we get that 
$$ f(t) + i(Hf)(t) = c_0 + 2\sum_{n=1}^{\infty} c_n e^{2\pi\,int} , $$
i.e. there is no negative spectrum anymore. That is exactly what Gabor needed, a transform that convert any real signal into one with the same real part, but without negative spectrum \cite{Ga}.
Surprisingly, this transform can be rewritten as an integral and as integral along the real line it becomes
$$ (Hf)(t) = v.p. \frac{1}{\pi}\int_{-\infty}^{\infty} \frac{f(\tau)}{t-\tau}\,d\tau . $$

The Hilbert transform is an important tool in optics. 
The optical implementation of the Hilbert transform was done in 1950 in two different ways. Kastler \cite{Kastler} used it for image processing, especially for edge enhancement. The mathematical reason for that is that the Hilbert transform is an singular integral operator and works like a derivative. Wolter \cite{Wolter} used it for spectroscopy which lead to further developments such as the fractional Fourier transform. All those applications make the Hilbert transform an important transform.
An higher dimensional analog of this useful transform are the Riesz transforms which can also be seen as a Hilbert transform.
In 1996 Lohmann \cite{LMZ} generalized the classical Hilbert transform to the fractional Hilbert transform, which can be implemented optically as a spatial filtering setup. 

It was widely believed, in the areas of optics, image analysis, and visual perception, that the Hilbert transform does not extend naturally and isotropically beyond one dimension. A transform that has these properties was in 2d introduced by Larkin \cite{Larkin} as a so-called spiral phase operator (in Fourier domain). It turns out that the spiral phase operator coincides with the monogenic signal.

Felsberg and Sommer \cite{FS} used the Riesz transforms or the higher dimensional Hilbert transform  to introduce the monogenic signal which is an higher dimensional analog of the analytic signal. The monogenic signal is different form the higher dimensional analytic signal defined in \cite{Ha}. The monogenic signal as well as the higher dimensional analytic signal have been proven to be useful in image processing (see for example \cite{BBRH}). 
In \cite{DNC} the fractional Hilbert transform was employed for image processing and specifically for edge detection.
In \cite{LMZ} the fractional Hilbert transform is deduced from an optical setup. 
$$ \tilde{H}_1(\nu) = \exp\left(+i\tfrac{\pi}{2}\right)S(\nu) + \exp\left(-i\tfrac{\pi}{2}\right)S(-\nu), $$
where $S(\nu)$ is a step function. The first fractional generalization is the filter:
$$ \tilde{H}_P(\nu) = \exp\left(+i\varphi\right)S(\nu) + \exp\left(-i\varphi\right)S(-\nu). $$
This can be rewritten as
$$ \tilde{H}_P(\nu) =\cos\varphi \tilde{H}_0(\nu) + \sin\varphi \tilde{H}_1(\nu). $$
For an input function $U(x)$ with a spatial spectrum $\tilde{U}(\nu),$ the output image is
\begin{align*}
V_P(x)  & = \int_{-\infty}^{\infty} \tilde{U}(\nu)\tilde{H}_P(\nu)\exp(2\pi i\nu x)\,d\nu  \\
        & = \cos\varphi \,V_0(x) + \sin\varphi\, V_1(x), 
\end{align*}
where $V_P$ is a mixture of the ordinary bright-field kind, $V_0,$ and the classical Hilbert transform $H_1.$
Mathematically, 
$$ S(\pm \nu) = \frac{1}{2}\left(1\pm\frac{\nu}{|\nu|}\right) $$
are projections related to the decomposition
$$ L^2(\mathbb{R}) = L^2(-\infty,0)\oplus L^2(0,\infty) $$
and the Hilbert transform is a convolution operator $H_1f= \frac{1}{\pi x}* f$ and the Fourier transform of the convolution kernel is $h(\nu) = \widehat{\frac{1}{\pi x}} = -i\, \mathrm{sgn\,}(\nu).$\\
We will generalizes these ideas into higher dimensions constructing a fractional Riesz-Hilbert transform in $\mathbb{R}^2.$
Later in \cite{Cu} three different generalizations of Gabor's signal were constructed, all of which reduce to the analytic signal when the angle is set to $\frac{\pi}{2}.$ It is of interest that this fractional Hilbert transform has the semigroup property, unlike the generalized Hilbert transform in \cite{Za}. 
Venkitaraman and Seelamantula considered several generalized concepts of analytic signals based on the fractional Hilbert transform \cite{VS}. 

There are to preliminary papers to that topic. In \cite{BE1} the construction of fractional Hilbert transform is done in Fourier domain based on rotations in $\mathbb{R}^3,$ which can be described by quaternions. The construction of the monogenic signal in \cite{BE2} is similar to that of the quaternionic monogenic signal in this paper. 

\subsection{Purpose of this paper}
In this paper we use a different approach to fractional Riesz-Hilbert transforms and fractional monogenic signals. A general construction of fractional transform can be done in the following way. Let $\mathcal{T}$ be transform such that $\mathcal{T}:\, H \to H$, where $H$ is a complex separable Hilbert with inner product $\langle .\, , .\rangle ,$ and if there is a complete set of orthonormal eigenfunctions $\phi_n$ with corresponding eigenvalues $\lambda_n,$ 
then any element in the space can be represented as
$$ f = \sum_{n=0}^{\infty} a_n\phi_n,\quad a_n = \langle f,\,\phi_n\rangle, \quad \mbox{so that}\quad  \mathcal{T}f = \sum_{n=0}^{\infty} a_n\lambda_n\phi_n . $$
The fractional transform can be defined as
$$ (\mathcal{T}^af)(\xi) = \sum_{n=0}^{\infty}a_n \lambda_n^a\phi_n(\xi) = \sum_{n=0}^{\infty} \lambda_n^a \langle f,\,\phi_n \rangle \phi_n(\xi) = \langle f,\, K_a(\xi,\, .) \rangle , $$
where 
$$ K_a(\xi, x) = \sum_{n=0}^{\infty} \overline{\lambda_n^a\phi_n(\xi)}\,\phi_n(x). $$
In such a case the transform should be compact operator but it is easily seen that the Hilbert transform as well as the Fourier transform are not compact operators in $L^2.$ But we can still use the approach because the Hilbert transform is involved in the Hardy-space decomposition of $L^2$ and actually the boundary values of monogenic functions in the upper/lower half-space are eigenfunctions of the Hilbert transform.

In optics and image processing the Hilbert transform is usually seen as the 1d Hilbert transform or the generalization to an analytic Hilbert transform in $\mathbb{C}^2$  \cite{Ha}.
Our paper is based on the Riesz transforms which are not considered separately but all together as a Hilbert transform in Clifford analysis. 
To distinguish between these two Hilbert transforms we will call our transform Riesz-Hilbert transform which is nothing else but the Hilbert transform in Clifford analysis.

The new  approach to construct fractional transforms is to substitute the eigenvalue $-1$ by $e^{-i\pi}$ (and $1 = e^0$) or $e^{-\underline{u}\pi},$ where $\underline{u}$ is a pure unit quaternion.  It will be proven that two (quaternionic) fractional Riesz-Hilbert transforms give rise to a fractional and quaternionic fractional monogenic signal respectively. These signals are modulated and rotated monogenic signals. It will be further proven that the (quaternionic) fractional Riesz-Hilbert transforms built a $C_0$-continuous semigroup.

The paper is organized as follows.  
After this introduction we provide the necessary facts about quaternionic analysis such as the algebra of quaternions itself, rotations in $\mathbb{R}^3$ and $\mathbb{R}^4$ and how they can be described by quaternions, the special case of isoclinic rotations. Further quaternionic analysis, the appropriate Dirac operator, the Hardy space decomposition and their projections. In section 3 we summarize the properties of the Riesz-Hilbert transform. The fractional Riesz-Hilbert transform and the fractional monogenic signal are constructed in section 4. In section 5 extend the theory to quaternionic fractional Riezs-Hilbert transforms and the quaternionic fractional monogenic signal.

\section{Preliminaries}
The Riesz-Hilbert transform can easily be defined in any Clifford algebra, but because we will use only 2D signals, the Hilbert transform will describe boundary values of monogenic functions in the upper half space $\mathbb{R}^3_+,$
and because we will use rotations in $\mathbb{R}^3,$ which are easily described by the aid of pure quaternions, we will use quaternions. 

The paper is organized as follows. After this introduction in Section 2 we real and complex quaternions, rotations in $\mathbb{R}^4$ and $\mathbb{R}^3,$ quaternionic analysis based on a Dirac operator and important properties of monogenic functions, Hardy spaces and the projections onto those spaces. In Section 3 we define the factional Hilbert transform and two special cases. The paper end with some concluding remarks.

\subsection{Quaternions}
\subsubsection{Real quaternions}
Let $\mathbb{H}$ denote the skew-field of quaternions with basis $\{\mathbf{1},\,\mathbf{i},\,\mathbf{j},\,\mathbf{k}\}$. An arbitrary element $q\in\mathbb{H}$ is given by
$$ q = q_0\mathbf{1}+q_1\mathbf{i}+q_2\mathbf{j}+q_3\mathbf{k}, $$
where $q_0,\,q_1,\,q_2,\,q_3\in\mathbb{R}$.
Addition and multiplication by a real scalar defined in component-wise fashion. The multiplication distributes over addition; $\mathbf{1}$ is the multiplication identity; and
\begin{align*}
	\mathbf{i}^2=\mathbf{j}^2=\mathbf{k}^2&=-\mathbf{1}, \\
\mathbf{i}\mathbf{j}=-\mathbf{j}\mathbf{i}=\mathbf{k},\quad \mathbf{j}\mathbf{k}=-\mathbf{k}\mathbf{j} & = \mathbf{i},\quad \mathbf{k}\mathbf{i}=-\mathbf{i}\mathbf{k}=\mathbf{j}.
\end{align*}
We decompose a quaternion into two parts that are called its scalar and vector parts. If
$ q = q_0\mathbf{1}+q_1\mathbf{i}+q_2\mathbf{j}+q_3\mathbf{k},$ we write
$$q =q_0 + \underline{q} = {\mathrm Sc\,}(q) + {\mathrm Vec\,}(q), $$
where $q_0 = q_0\mathbf{1}$ is the scalar part and $\underline{q}=q_1\mathbf{i}+q_2\mathbf{j}+q_3\mathbf{k} $ the vecor part of $q$.
On $\mathbb{H}$ a conjugation is defined as $\overline{q}= q_0\mathbf{1}-\underline{q}$ for all $q\in\mathbb{H}.$ The length or norm  $|q|$ of $q\in \mathbb{H}$ is then
$$ |q| = \sqrt{q_0^2+q_1^2+q_2^2+q_3^2} = \sqrt{q\overline{q}} = \sqrt{\overline{q}q}. $$
All $q\in\mathbb{H}\backslash\{0\}$ are invertible with inverse $$q^{-1}=\frac{\overline{q}}{|q|}.$$


\subsection{Rotations}
We identify $\mathbb{H}$ with the four dimensional Euclidean space $\mathbb{R}^4$ by associating $q\in\mathbb{H}$ with the vector $(q_0,\,q_1,\,q_2,\,q_3)\in\mathbb{R}^4$.

If $q\in\mathbb{H}$ has $|q|=1,$ then we call $q$ a unit quaternion. If $\underline{u}$ is a pure unit quaternion then $\underline{u}^{-1}=-\underline{u}.$ Two pure unit quaternions $\underline{u}$ and $\underline{v}$ are orthogonal if and only if 
$$ \underline{u}\underline{v}+\underline{v}\underline{u}=0. $$
If $q$ is a unit quaternion there is a real number $\varphi$ and a pure unit quaternion $\underline{u}$ such that 
$$ q = \mathbf{1}\cos\varphi + \underline{u}\sin\varphi .$$
Since $\underline{u}^2=-1,$ the power series expansion leads to
\begin{align*}
e^{\underline{u}\varphi}  & = \sum_{n=0}^{\infty} \frac{(\underline{u}\varphi)^n}{n!}  = \sum_{l=0}^{\infty} \frac{(\underline{u}\varphi)^{2l}}{(2l)!} + \sum_{l=0}^{\infty} \frac{(\underline{u}\varphi)^{2l+1}}{(2l+1)!} \\
& = \sum_{l=0}^{\infty} \frac{(-1)^l \varphi^{2l}}{(2l)!} + \underline{u}\sum_{l=0}^{\infty} \frac{(-1)^l\varphi^{2l+1}}{(2l+1)!} 
= \mathbf{1}\cos\varphi + \underline{u}\sin\varphi , 
\end{align*}
providing equivalent representations for a unit quaternion
$$ q = q_0 + \underline{q} = \mathbf{1}\cos\varphi + \underline{u}\sin\varphi = e^{\underline{u}\varphi}. $$
It should be mentioned that neither $\underline{u}$ nor $\varphi$ is uniquely determined. When $q\not=\pm\mathbf{1}$ then $\sin\varphi = \pm |\underline{q}|$ and $\underline{u}=\pm\tfrac{\underline{q}}{|\underline{q}|}$. When $q = \pm \mathbf{1},$ $\underline{u}$ can be any pure unit quaternion. 

\subsubsection{Rotations in $\mathbb{R}^3$}
Let $q$ be a quaternion. Then $L(q):\,\mathbb{H}\to\mathbb{H}$ and $R(q):\,\mathbb{H}\to\mathbb{H}$ are defined as follows
$$ [L(q)](x)=qx,\quad [R(q)](x)=xq,\quad x\in\mathbb{H}. $$
If $q$ is a unit quaternion, then both $L_q$ and $R_q$ are orthogonal transformations of $\mathbb{H}.$  Thus, for unit quaternions $p$ and $q$, the mapping $C_{p,q}:\,\mathbb{H}\to\mathbb{H}$ defined by
$$ C_{p,q} = L(p) \circ R(q) = R(q) \circ L(p)  $$
is also an orthogonal transform of $\mathbb{H}$ and $C_{p_1,q_1}\circ C_{p_2,q_2}=C_{p_1p_2,q_2q_1}.$

\begin{thm}[\cite{WW}] If $q$ is a unit quaternion, then there exist a pure unit quaternion $\underline{u}$ and a real scalar $\varphi$ such that $q = e^{\underline{u}\varphi}.$ The transform $C:\,\mathbb{R}^3\to\mathbb{R}^3$ defined by $C(x)=qx\overline{q}$ is a rotation in the plane orthogonal to $\underline{u}$ through an angle $2\varphi.$
\end{thm}
	
\subsubsection{Rotations in $\mathbb{R}^4$}

\begin{thm}[\cite{WW}] Let $p=e^{\underline{u}\varphi}$ and $q=e^{\underline{v}\psi},$ where $\underline{u}$ and $\underline{v}$ are pure unit quaternions. The orthogonal transform $C_{p,q}$ of $\mathbb{H}$ is a product of two rotations in orthogonal planes. 
\end{thm}


This theorem is presumably due to Ball; in \cite{Ball} the author does not mention it explicitly as a theorem, but nevertheless gives a proof. However,
Ball's proof is slightly incomplete, a complete proof is given by Mebius  in \cite{Mebi}.
Hence, the multiplication with an unit quaternion from the right (or left) only describes a (right- or left-) isoclinic rotation in $\mathbb{R}^4.$

We will use a multiplication from the right-hand side. 

Let $P$
be an arbitrary 4D point, represented as a quaternion
$P=w\mathbf{1}+x\mathbf{i}+y\mathbf{j}+z\mathbf{k}$. Let $p$ and $q$
be unit quaternions. Consider the left- and right-multiplication mappings $P\mapsto pP $ and $P \mapsto Pq.$
Both mappings have the property of rotating all half-lines originating from $\mathcal{O}$ through the same angle ($\arccos p_0$ and $\arccos q_0$ respectively); such rotations are denoted as $\mathbf{isoclinic}$ \cite{WW}.

\begin{defn}[Isoclinic rotation] In a 4-dimensional space a rotation is said isoclinic if all its 2 angles are equal (up to the sign).
\end{defn}

Because the left- and the right multiplication are different from each other resulting in different rotations we have to distinguish between left- and right-isoclinic rotations.
Conversely, an isoclinic 4D rotation about $\mathcal{O}$ (different from the non-rotation $I$ and from the central reversion $-I$) is represented by either a left-multiplication or a right-multiplication by a unique unit quaternion. Hence, the multiplication with an unit quaternion from the right (or left) only describes a right- (or left-) isoclinic rotation in $\mathbb{R}^4.$

\subsection{Function spaces}
We consider functions in the Lebesgue spaces $L^p.$ A real-,  complex-, or quaternionic-valued function $f$ belongs to $L^p(\mathbb{R}^n,\,\mathbb{R}),\, L^p(\mathbb{R}^n,\,\mathbb{C})$  or $L^p(\mathbb{R}^n,\,\mathbb{H})$ if
$$ \left(\int_{\mathbb{R}^n} |f(\underline{x})|^p\,d\underline{x}\right)^{\frac{1}{p}},\ 1<p<\infty . $$
The space $L^2(\mathbb{R}^n,\mathbb{H})$ is a Hilbert space with inner product
$$ \langle f,\, g \rangle = {\rm Sc\,} \int \overline{f(\underline{x})}^{\mathbb{CH}}\,g(\underline{x})\\,d\underline{x} =  \int \sum_{j=0}^3 \overline{f_j(\underline{x})}^{\mathbb{C}}\,g_j(\underline{x})\,e_j^2\,d\underline{x} , $$
where $ \overline{\ }^{\mathbb{CH}} $ denotes the complex and quaternionic conjugation and $\overline{\ }^{\mathbb{C}}$ the complex conjugation.


\subsection{Quaternionic analysis}
The quaternionic analysis presented in this section is based on the nice presentation of Clifford analysis in \cite{De}.
\subsubsection{Dirac operator}
The Dirac operator is defined as the first order  linear differential operator
$$ D_{x} = \mathbf{1}\partial_{x_0}+\mathbf{i}\partial_{x_1}+\mathbf{j}\partial_{x_2}+\mathbf{k}\partial_{x_3}. $$

\begin{defn}[Monogenic functions] Let $\Omega \subseteq \mathbb{R}^{4}$ be open and let $f$ be a $C^1$-function in $\Omega$ which is (real or complex ) quaternionic-valued. Then $f$ is {\it left monogenic} or  ({\it right monogenic}) in $\Omega$ if in $\Omega$ 
$$ D_xf = 0 \quad \mbox{or} \quad fD_x = 0.$$
\end{defn}
The connection between monogenic and harmonic functions is due to the fact that $\Delta_x = -\partial_x^2=\overline{D}_xD_x.$

\subsubsection{Integral formulae}
Let $\Omega\subset \mathbb{R}^{4}$ be open.
Because the Dirac operator $D_x$ is a first order linear differential operator with constant coefficients there exists a fundamental solution. 
\begin{lem}\label{fundamental solution}
A fundamental solution is given by
$$ E(x) = \frac{1}{2\pi^2} \frac{\overline{x}}{|x|^{4}} = \frac{1}{2\pi^2} \frac{(x_0-x_1\mathbf{i}-x_2\mathbf{j}-x_3\mathbf{k})}{|x|^{4}}. $$ 

$E$ has the following properties.
\begin{enumerate} 
 \item $E$ is $\mathbb{H}\sim\mathbb{R}^{4}$-valued and belongs to $L^1_{loc}(\mathbb{R}^{4}).$
 \item $E$ is left and right monogenic in $\mathbb{R}^{4}\backslash\{0\}$ 
\item and $\lim\limits_{|x|\to\infty} E(x) = 0.$
\item $D_xE = ED_x = \delta(x),$ $\delta(x)$ being the classical $\delta$-function in $\mathbb{R}^{4}.$
\end{enumerate}
\end{lem}

\subsubsection{Hardy spaces}

\begin{defn}[Hardy spaces] Let $1<p<\infty.$ Then 
$$ H^p(\mathbb{R}^3_{\pm}) = \{ F \ \mbox{monogenic in\ } \mathbb{R}^4_{\pm}: \sup_{\varepsilon > 0} \int_{\mathbb{R}^3} |F(y\pm \varepsilon)|^p\,dy<\infty\} $$
is the Hardy space of (left) monogenic functions in $\mathbb{R}^3_{\pm}.$ 
\end{defn}


\begin{defn}[Integral operators] For $f\in L^p(\mathbb{R}^3)$ and $x\in\mathbb{R}^{4}\backslash \mathbb{R}^3 ,$ 
$$ \mathcal{C}f(x) = \int_{\mathbb{R}^3} E(x-y)f(y)\,dy $$
is the Cauchy transform of $f.$ \\
For $f\in L^p(\mathbb{R}^4)$ and a.e. $x\in\mathbb{R}^3 ,$
\begin{align*}
(Hf)(x)& = 2\, p.v. \int_{\mathbb{R}^3} E(x-y)f(y)\,dS(y) \\
&  = 2 \lim_{\varepsilon \to 0+} \int_{y\in\mathbb{R}^3 : |x-y|>\varepsilon} E(x-y)f(y)\,dy
\end{align*}
is the Hilbert transform of $f$.
\end{defn}
the Riesz-Hilbert transform is a convolution operator, i.e. $Hf = p.v. 2E(\underline{x})\ast f.$ We call the Fourier transform of the convolution kernel Fourier symbol of the operator. Hence the Riesz-Hilbert transform has Fourier symbol $\frac{i\underline{\xi}}{|\underline{\xi}|}.$

\begin{thm}[\cite{McI}] Let $f\in L^p(\mathbb{R}^3),\ 1<p<\infty.$ Then
\begin{enumerate}
\item $\mathcal{C}f \in H^p(\mathbb{R}^3_{\pm}).$
\item $\mathcal{C}f$ has non-tangential limits $(\mathcal{C}^{\pm})$ at almost all $x^*\in\mathbb{R}^3.$
\item Putting
$$ \mathbb{P}^+f(x^*) = (\mathcal{C}^+f)(x^*) \quad \mbox{and}\quad \mathbb{P}^-f(x^*) = -(\mathcal{C}^-f)(x^*)$$
then $\mathbb{P}^{\pm}$ are bounded projections in $L^p(\mathbb{R}^3).$
\item (Plemelj-Sokhotzki formulae). For a.e. $x^*\in\mathbb{R}^3,$
$$ \mathbb{P}^+f(x^*) = \tfrac{1}{2}\left(f(x^*)+(Hf)(x^*)\right) \quad \mbox{and}\quad \mathbb{P}^-f(x^*) =  \tfrac{1}{2}\left(f(x^*)- (Hf)(x^*)\right)$$
whence 
$$ 1 = \mathbb{P}^+ + \mathbb{P}^- \quad \mbox{and}\quad  H = \mathbb{P}^+ - \mathbb{P}^- .$$
In particular $H$ is a bounded linear operator on $L^p(\mathbb{R}^3)$ and, putting $L^{p,\pm}(\mathbb{R}^3) = \mathbb{P}^{\pm}L^p(\mathbb{R}^3),$ leads to the decomposition into Hardy spaces
$$ L^p(\mathbb{R}^3) =  L^{p,+}(\mathbb{R}^3)  \oplus  L^{p,-}(\mathbb{R}^3) .$$
\end{enumerate}
\end{thm}

Using that 
$$ \mathcal{F}\left(\frac{2x_j}{2\pi^2|\underline{x}|^{4}}\right) = -i\frac{\xi_j}{|\underline{\xi}|}, \quad j=1,\ldots,3, $$
and set \cite{McI}
$$ \chi_{\pm}(\underline{\xi}) = \tfrac{1}{2}\left( 1 \pm i \frac{\underline{\xi}}{|\underline{\xi}|}\right),$$
it should be noticed that 
$$ \chi_{\pm}^2(\underline{\xi}) = \chi_{\pm}(\underline{\xi}) \quad \mbox{and} \quad \chi_+(\underline{\xi})+\chi_- (\underline{\xi})= 1,\quad \chi_+(\underline{\xi})\chi_-(\underline{\xi}) = \chi_-(\underline{\xi})\chi_+(\underline{\xi}) = 0.$$
Which means that $\chi_{\pm}$ are projections and zero divisors. 
That is equivalent to a decomposition in Fourier space:
$$ \hat{f} = \frac{1}{2}\left(1+i\frac{\underline{\xi}}{|\underline{\xi}|}\right)\hat{f}+ \frac{1}{2}\left(1-i\frac{\underline{\xi}}{|\underline{\xi}|}\right)\hat{f}.$$

The boundary values of monogenic function in upper half-space are characterized in the next theorem.

\begin{thm} For $f\in L^p(\mathbb{R}^3)$ the following statements are equivalent:
\begin{enumerate}
\item The non-tangential limit of $\mathcal{C}f$ is a.e. equal to $f,$
\item $Hf = f,$
\item $\mathcal{F}f = \chi_+\,\mathcal{F}f.$
\end{enumerate}
and characterizes boundary values of monogenic functions.
\end{thm}

\section{Properties of the Riesz-Hilbert transforms}
The Riesz-Hilbert transform  $H$ is a convolution-type operator and a linear combination of Riesz transforms $R_j:$
$Hf(\underline{x}) = \sum_{j=1}^3R_jf(\underline{x})e_j$ and the Fourier transform is $\mathcal{F}(Hf) = \frac{i\underline{\xi}}{|\underline{\xi}|}\hat{f},$ where $\hat{f}$ denotes the Fourier transform of $f$. The Riesz-Hilbert transform has the following properties:

\begin{enumerate}
\item \textit{Shift-invariance}: $\mathcal{S}_{\tau}(Hf)(\underline{x}) = H(\mathcal{S}_{\tau}f)(\underline{x}),\ \tau \in \mathbb{R}^3,$
\item \textit{Scale-invariance}: $\mathcal{D}_{\sigma}(Hf)(\underline{x}) = H(\mathcal{D}_{\sigma}f)(\underline{x}),\ \sigma \in \mathbb{R}^+,$
\item \textit{Self-Reversibility}: $\forall f \in L^2(\mathbb{R}^3,\mathbb{H}): H^2f(\underline{x}) =-\sum\limits_{i=1}^3 (R^2_jf)(\underline{x}) = f(\underline{x}).$
\item \textit{Energy preservation:} For all $f,\,g\in L^2(\mathbb{R}^3,\,\mathbb{H}): \langle Hf,\, Hg\rangle = \langle f,\, g\rangle .$
\item \textit{Orthogonality I}: The real-valued function $f\in L^2(\mathbb{R}^3,\,\mathbb{R})$ is orthogonal to its the Riesz-Hilbert transform $Hf\in L^2(\mathbb{R}^3,\mathbb{H}),$ i.e. $ \langle f,\, Hf \rangle = 0.$ 
\item \textit{Orthogonality II}: If $f,\,g\in L^2(\mathbb{R}^3, \,\mathbb{R}) $ are such that $\langle f,\, g \rangle = 0,$ then $\langle Hf,\,Hg\rangle = 0.$
\end{enumerate}
Proof: Because the Riesz transforms are linear convolution-type integrals the Riesz-Hilbert transform is linear. The Riesz-Hilbert transform is invariant under actions of the Vahlen group $V(n)$ \cite{Ryan}. Specific Vahlen matrices are 
\begin{eqnarray*}
\left(\begin{matrix}0 &  -1 \\ 1 & 0 \end{matrix}\right),\quad \left(\begin{matrix} a & 0 \\ 0 & 
{\tilde{a}}^{-1} \end{matrix}\right), \mbox{  where   } a\in Pin(n),\hspace*{2cm} \\
\left(\begin{matrix} 1 & v \\ 0 & 1 \end{matrix}\right), \mbox{  where  } v \in \mathbb{R}^n, \mbox{  and  } \left(\begin{matrix} \lambda & 0 \\ 0 & \lambda^{-1} \end{matrix}\right), \mbox{  where  } \lambda \in \mathbb{R}^+.
\end{eqnarray*}
These matrices are generators of the group $V(n)$ \cite{V}. The pin group $Pin(n)$ is the group of all even products of pure vectors or pure quaternions. The matrices correspond respectively to Kelvin inversion, orthogonal transformations, translation and dilation.
It is easily seen in Fourier domain that $\mathcal{F}(H^2f) = \left(\frac{i\underline{\xi}}{|\underline{\xi}|}\right)^2 \hat{f} = \hat{f}. $ The orthogonality of the real-valued signal $f$ to its Riesz-Hilbert transform is obvious from $\mathrm{Sc\,} (f\, (Hf)) = 0$. 
Because of 
$$ \langle Hf,\, Hg \rangle = -\sum_{j=1}^3 \int_{\mathbb{R}^3} (R_jf)(\underline{x})(R_jg)(\underline{x})\,d\underline{x} = -\sum_{j=1}^3\int_{\mathbb{R}^3}  \frac{i\xi_j}{|\underline{\xi}|}\hat{f}(\underline{\xi}) \frac{i\xi_j}{|\underline{\xi}|}\hat{g}(\underline{\xi}) \,d\underline{\xi} = \langle f,\, g \rangle $$
the orthogonality of $f$ and $g$ implies the orthogonality of their Riesz-Hilbert transforms. 
\begin{lem} For $f,\,g\in L^2(\mathbb{R}^3,\,\mathbb{H})$ we have 
$\langle Hf,\, g\rangle = \langle f,\,Hg \rangle .$
\end{lem}
Proof: We have
\begin{eqnarray*}  \int_{\mathbb{R}^3} \overline{\widehat{Hf}}^{\mathbb{CH}}\,\hat{g}\, d\underline{\xi} = \int_{\mathbb{R}^3} \overline{\hat{f}(\xi)}^{\mathbb{CH}}\,\frac{i\underline{\xi}}{|\underline{\xi}|}\,\hat{g}(\underline{\xi})\,d\underline{\xi} = \int_{\mathbb{R}^3} \overline{\hat{f}}^{\mathbb{CH}}\widehat{Hg}\,d\underline{\xi}
\end{eqnarray*}
and thus $ \langle Hf,\,g\rangle =  \langle f,\, Hg \rangle .$
 



\section{Fractional Riesz-Hilbert transform}
It is a general question how to construct a fractional transformation. The answer could be given in terms of a eigenvalue or singular value decompostion. Unfortunately, these theories work very well for compact operators only. The usual transforms are not compact operators. Nevertheless, we will use the basic idea of eigenvalue decomposition to obtain a definition of fractional transform. \\[1ex]
The idea by Lohmann to construct the fractional Hilbert transform leads to a norm continuous $C_0$-semigroup of Hilbert operators. The construction can be interpreted as rotations in Fourier domain. 
Because $H^2=I$ there are only 2 eigenvalues $\pm 1$. The eigenspace for $\lambda = 1$ is the Hardy space $L^{2,+}(\mathbb{R}^3)$ and for $\lambda =-1$ it is the Hardy space $L^{2,-}(\mathbb{R}^n).$

We have 
$$ \varphi = \frac{1}{2}\left(I+H\right)\varphi + \frac{1}{2}\left(I-H\right)\varphi $$
and
$$ H\varphi = \frac{1}{2}\left(I+H\right)\varphi + (-1) \frac{1}{2}\left(I-H\right)\varphi .$$
Therefore we define
\begin{align*}
H^{\alpha}\varphi & = \frac{1}{2}\left(I+H\right)\varphi + e^{-i\pi\alpha} \frac{1}{2}\left(I-H\right)\varphi  \\
                             & = e^{-i\frac{\pi}{2}\alpha}\frac{1}{2}\left(e^{i\frac{\pi}{2}\alpha}\left(I+H\right)\varphi + e^{-i\frac{\pi}{2}\alpha} \left(I-H\right)\right)\varphi \\
                            & = e^{-i\frac{\pi}{2}\alpha}\left(\cos\left(\frac{\pi}{2}\alpha\right) I + i \sin\left(\frac{\pi}{2}\alpha\right) H \right)\varphi \\
 \iff i^{\alpha}H^{\alpha}\varphi = (iH)^{\alpha}\varphi & = \cos\left(\frac{\pi}{2}\alpha\right) I + i \sin\left(\frac{\pi}{2}\alpha\right) H\varphi.
\end{align*}
Here, we identify $e^{-i\tfrac{\pi}{2}\alpha}$ with $i^{\alpha}.$ From that we obtain the following definition.
\begin{defn}[Fractional Riesz-Hilbert transform] 
The fractional Riesz-Hilbert transform is defined as
\begin{align}
 & H^{\alpha} := e^{-i\tfrac{\pi}{2}\alpha}\left(\cos(\tfrac{\pi}{2}\alpha)I + i\sin(\tfrac{\pi}{2})H)\right)   \\
\mbox{and }\quad   & \mathcal{H}^{\alpha} = (iH)^{\alpha}  := \cos(\tfrac{\pi}{2}\alpha)I + i\sin(\tfrac{\pi}{2})H . \label{FHRT}
\end{align}
\end{defn}
\begin{thm}
The fractional Riesz-Hilbert transforms $H^{\alpha},\,\mathcal{H}^{\alpha}: L^p(\mathbb{R}^3,\,\mathbb{H}) \to L^p(\mathbb{R}^3,\,\mathbb{H}),\\ 1<p<\infty, \,\alpha\in\mathbb{R},$ are linear, shift and scale invariant and fulfill the following properties
\begin{enumerate}
\item the inverse of $H^{\alpha}$ is $H^{-\alpha}$ and the inverse of $\mathcal{H}^{\alpha}$ is $\mathcal{H}^{-\alpha}, \alpha \in \mathbb{R},$ 
\item $H^{\alpha}$ is $2$-periodic in $\alpha,$ i.e. $H^{\alpha +2} = H^{\alpha},$ whereas $\mathcal{H}^{\alpha}$ is $4$-periodic in $\alpha,$ i.e. $\mathcal{H}^{\alpha + 4} = \mathcal{H}^{\alpha},\ \alpha\in\mathbb{R}.,$
\item If $f,\,g \in L^2(\mathbb{R}^3,\,\mathbb{R})$ such that $\langle f,\,g \rangle = 0$ then $$\langle H^{\alpha}f,\, H^{\alpha}g\rangle = \langle \mathcal{H}^{\alpha}f,\, \mathcal{H}^{\alpha}g\rangle = 0.$$
\end{enumerate}

\end{thm}
Because $H$ maps $L^p(\mathbb{R}^3, \mathbb{H})$ onto itself both operators are linear and bounded operators in $L^p(\mathbb{R}^3,\mathbb{H}).$ 
The Hilbert operator as a singular integral operator is a linear bounded operator in $\mathcal{L}(L^p(\mathbb{R}^n)),\,1<p<\infty .$ 
Further, we have
$$ \mathcal{H}^{\alpha} = \cos\left(\frac{\pi}{2}\alpha\right)I + i \sin\left(\frac{\pi}{2}\alpha\right) H = \cos\left(\frac{\pi}{2}(\alpha+4)\right)I + i \sin\left(\frac{\pi}{2}(\alpha+4)\right) H = \mathcal{H}^{\alpha +4}, $$
because sine and cosine function are $2\pi$-periodic, and
\begin{align*} 
H^{\alpha} & = e^{-i\frac{\pi}{2}\alpha}\left(\cos\left(\frac{\pi}{2}\alpha\right) I + i \sin\left(\frac{\pi}{2}\alpha\right) H \right)  = -e^{-i\frac{\pi}{2}\alpha}\left(-\cos\left(\frac{\pi}{2}\alpha\right) I - i \sin\left(\frac{\pi}{2}\alpha\right) H \right) \\
                  & = e^{-i\frac{\pi}{2}(\alpha+2)}\left(\cos\left(\frac{\pi}{2}(\alpha+2)\right) I + i \sin\left(\frac{\pi}{2}(\alpha+2)\right) H \right) = H^{\alpha +2}.
\end{align*}

\begin{rem} Some examples: $\mathcal{H}^0 = H^0 = I,$ $ \mathcal{H}^1 = iH $ and $H^1 = H,$ $\mathcal{H}^2= -I,$ where as $H^2 = I,$ and $\mathcal{H}^3 = -iH,$ $H^3 = H,$ $\mathcal{H}^4 = I.$
\end{rem}
\begin{rem} The operator $\mathcal{H}^{\alpha} = e^{i\tfrac{\pi}{2}\alpha}\tfrac{1}{2}(I+H) + e^{-i\tfrac{\pi}{2}\alpha}\tfrac{1}{2}(I-H)$ is a convolution operator with Fourier symbol
$$  e^{i\tfrac{\pi}{2}\alpha}\tfrac{1}{2}\left(1+\tfrac{i\underline{\xi}}{|\underline{\xi}|}\right) + e^{-i\tfrac{\pi}{2}\alpha}\tfrac{1}{2}\left(1-\tfrac{i\underline{\xi}}{|\underline{\xi}|}\right) = e^{i\tfrac{\pi}{2}\alpha}\chi_+(\underline{\xi}) +  e^{-i\tfrac{\pi}{2}\alpha}\chi_-(\underline{\xi}).$$
That means that the construction of $\mathcal{H}^{\alpha}$ is done in a similar way as the construction of the fractional Hilbert transform by Lohmann. $\chi_{\pm}$ are the step functions in higher dimensions.
\end{rem}
\begin{rem} In Fourier domain we get another nice construction of the Fourier symbol of the fractional Riesz-Hilbert transform. The Fourier symbol of $\mathcal{H} = iH $ is 
$$i\frac{i\underline{\xi}}{|\underline{\xi}|} = -\frac{\underline{\xi}}{|\underline{\xi}|} = \cos\left(\tfrac{\pi}{2}\right) - \frac{\underline{\xi}}{|\underline{\xi}|} \sin\left(\tfrac{\pi}{2}\right) = e^{-\frac{\underline{\xi}}{|\underline{\xi}|}\frac{\pi}{2}}$$
and hence the Fourier symbol of $\mathcal{H}^{\alpha}$ is
$$ \left(i\frac{i\underline{\xi}}{|\underline{\xi}|}\right)^{\alpha} = \left(-\frac{\underline{\xi}}{|\underline{\xi}|}\right)^{\alpha} = e^{-\frac{\underline{\xi}}{|\underline{\xi}|}\frac{\pi}{2}\alpha} = \cos\left(\tfrac{\pi}{2}\alpha\right) + i \,\frac{i\underline{\xi}}{|\underline{\xi}|} \sin\left(\tfrac{\pi}{2}\alpha\right). $$
It has to be mentioned that this is more a definition than a property because a fractional power of a quaternion is not uniquely defined. For example the equation $\underline{u}^2 = -1$ has infinite many solutions.
\end{rem}
The considerations on the Fourier symbol implies that the fractional Riesz-Hilbert transform should have a semigroup property.
\begin{defn}[Semigroup] A family $T=\{T_t\}_{t\geq 0}$ of bounded linear operators acting on a Banach space $E$ is called a $C_0$-semigroup if the following three properties are satisfied:
\begin{enumerate}
 \item $T_0 = I,$
 \item $T_{s+t} = T_sT_t$ for all $s,\,t\geq 0,$
 \item $\lim\limits_{t\to 0} T_tx = x $ for all $x\in E.$
\end{enumerate}
If the stronger condition $ \lim\limits_{t\to 0} ||T_t - T||_E = 0$ is satisfied the group is called norm continuous.
\end{defn}
\begin{thm} The Riesz-Hilbert transform generates a norm continuous $C_0$-semigroup in $L^p(\mathbb{R}^3,\mathbb{H}).$ We have
$$ e^{i\frac{\pi}{2}\alpha\, H} =\cos\left(\frac{\pi}{2}\alpha\right)I + i \sin\left(\frac{\pi}{2}\alpha\right) H = \mathcal{H}^{\alpha}. $$
The operator $I-H$ also generates a norm continuous $C_0$-semigroup in $L^p(\mathbb{R}^3,\mathbb{H}).$ We have
$$ e^{-i\pi\alpha \frac{1}{2}(I-H)} 
     = e^{-i\frac{\pi}{2}\alpha}\left(\cos\left(\frac{\pi}{2}\alpha\right) I + i \sin\left(\frac{\pi}{2}\alpha\right) H \right)  = H^{\alpha}. $$
\end{thm}
Proof: We have
\begin{align*}
e^{i\frac{\pi}{2}\alpha\, H} & = \sum_{l=0}^{\infty} i^{2l}\left(\frac{\pi}{2}\alpha\right)^{2l}\frac{1}{(2l)!} H^{2l} + i^{2l+1}\left(\frac{\pi}{2}\alpha\right)^{2l+1}\frac{1}{(2l+1)!} H^{2l+1} \\
 & = \sum_{l=0}^{\infty} (-1)^l\left(\frac{\pi}{2}\alpha\right)^{2l}\frac{1}{(2l)!} I + i (-1)^l\left(\frac{\pi}{2}\alpha\right)^{2l+1}\frac{1}{(2l+1)!} H \\
 & = \cos\left(\frac{\pi}{2}\alpha\right)I + i \sin\left(\frac{\pi}{2}\alpha\right) H .
\end{align*}
and 
\begin{align*}
e^{-i\pi\alpha \frac{1}{2}(I-H)} & = e^{-i\frac{\pi}{2}\alpha I} e^{i\frac{\pi}{2}H}  
     = e^{-i\frac{\pi}{2}\alpha}\left(\cos\left(\frac{\pi}{2}\alpha\right) I + i \sin\left(\frac{\pi}{2}\alpha\right) H \right)
    \end{align*}
Furthermore,
\begin{align*}
\lim_{\alpha \to 0} || \mathcal{H}^{\alpha}f - f||_{L^p} & = \lim_{\alpha\to 0} || \cos\left(\tfrac{\pi}{2}\alpha\right)If + i \sin\left(\tfrac{\pi}{2}\alpha\right) Hf - f ||_{L^p} \\
 & \leq \lim_{\alpha\to 0} \left(|\cos(\tfrac{\pi}{2}\alpha) - 1| \,||f||_{L^p} + |\sin(\tfrac{\pi}{2}\alpha)|\, ||Hf||_{L^p}\right) = 0 \\
\mbox{and  } \lim_{\alpha \to 0} || H^{\alpha}f - f||_{L^p} & = \lim_{\alpha \to 0} || e^{-i\tfrac{\pi}{2}\alpha}\mathcal{H}^{\alpha}f - f||_{L^p}  \\
& = \lim_{\alpha \to 0} || (e^{-i\tfrac{\pi}{2}\alpha}-1)\mathcal{H}^{\alpha}f + \mathcal{H}^{\alpha}f- f||_{L^p} \\
& = \lim_{\alpha \to 0} || (\cos({\tfrac{\pi}{2}\alpha})-1)\mathcal{H}^{\alpha}f - i \sin (\tfrac{\pi}{2}\alpha) \mathcal{H}^{\alpha}f+ \mathcal{H}^{\alpha}f- f||_{L^p} \\
& \hspace*{-2cm}=\lim_{\alpha \to 0} \left(|\cos({\tfrac{\pi}{2}\alpha})-1|\,||\mathcal{H}^{\alpha}f||_{L^p} + |\sin (\tfrac{\pi}{2}\alpha)|\,||\mathcal{H}^{\alpha}f||_{L^p}+ ||\mathcal{H}^{\alpha}f- f||_{L^p} \right) = 0.
\end{align*}

\begin{defn}[Fractional monogenic signal] For $f\in L^p(\mathbb{R}^3,\mathbb{R})$ the fractional monogenic signal is defined as
$$ M^{\alpha}f := (I-e^{i\frac{\pi}{2}\alpha}\mathcal{H}^{\alpha})f = (I+e^{i\pi(1-\alpha)}H^{\alpha})f . $$
\end{defn}


\begin{thm}[Properties of the fractional monogenic signal] 
Let $f\in L^p(\mathbb{R}^3,\mathbb{R}),$ then the fractional monogenic signal is linear, shift- and scale-invariant and possesses the following properties:
\begin{enumerate}
\item $M^{\alpha}f \in H^+(\mathbb{R}^3),$ i.e. boundary values of a monogenic function in the upper half space.
\item The fractional monogenic signal is a rotated and modulated version of the monogenic signal.
\item $||M^{\alpha}f(\underline{x})||_{L^p} = |\sin(\frac{\pi}{2}\alpha)|\, ||Mf||_{L^p}.$ 
\item The fractional monogenic signal is not orthogonal to the signal and
$$ Mf = (I+H)f = (1+i\cot (\tfrac{\pi}{2}\alpha))f - i \csc (\tfrac{\pi}{2}\alpha) H^{\alpha}f. $$
\item The fractional monogenic signal is $2$-periodic.
\end{enumerate}
\end{thm}
Proof: We have
\begin{align*}
M^{\alpha}f & = (I - e^{i\tfrac{\pi}{2}\alpha} \mathcal{H}^{\alpha})f = (I+e^{i\pi(1-\alpha)}H^{\alpha})f \\
& = I - \left(\cos(\tfrac{\pi}{2}\alpha) + i\sin(\tfrac{\pi}{2}\alpha)\right)\left(\cos(\tfrac{\pi}{2}\alpha) + i\sin(\tfrac{\pi}{2}\alpha)H\right) \\
& = \sin^2(\tfrac{\pi}{2}\alpha) - i\sin(\tfrac{\pi}{2}\alpha)\cos(\tfrac{\pi}{2}\alpha) +\left(-i\sin(\tfrac{\pi}{2}\alpha)\cos(\tfrac{\pi}{2}\alpha) + \sin^2(\tfrac{\pi}{2}\alpha)\right) H \\
&  = -2i \sin(\tfrac{\pi}{2}\alpha)e^{i\tfrac{\pi}{2}\alpha}\tfrac{1}{2}(I+H) = -i\sin(\tfrac{\pi}{2}\alpha)e^{i\tfrac{\pi}{2}\alpha}Mf,
\end{align*}
i.e.  $M^{\alpha}f$ is up to a constant depending on $\alpha$ equal to the monogenic signal. We have 
$$ |M^{\alpha}f(\underline{x})| = |\sin(\frac{\pi}{2}\alpha)|\,|Mf(\underline{x})| \quad \mbox{and}\quad ||M^{\alpha}f||_{L^p} = |\sin(\frac{\pi}{2}\alpha)|\, ||Mf||_{L^p} . $$
Because the fractional Riesz-Hilbert transform can be expressed by the Riesz-Hilbert transform we get
\begin{align*}
Mf & = (I+H)f  
= \left(f - \left(e^{i\tfrac{\pi}{2}\alpha} \mathcal{H}^{\alpha}f - \cos(\alpha\tfrac{\pi}{2})f\right) \frac{i}{\sin(\tfrac{\pi}{2}\alpha)}\right) \\
& = \left(1+ i\cot(\tfrac{\pi}{2}\alpha) - i H^{\alpha}\csc(\tfrac{\pi}{2})\right)f. 
\end{align*}
We have
$$ M^{\alpha}f = -i\sin(\tfrac{\pi}{2}\alpha)e^{i\tfrac{\pi}{2}\alpha}Mf = -i\sin(\tfrac{\pi}{2}(\alpha+2))e^{i\tfrac{\pi}{2}(\alpha+2)}Mf = M^{\alpha + 2}f.$$
Therefore, the fractional monogenic is $2$-periodic.

\section{Quaternionic fractional Riesz-Hilbert transform}

To define the fractional Riesz-Hilbert transform we used  $-1 = e^{-i\pi}$. \\ But we have also $-1 = e^{-\underline{u}\pi} = \cos(\pi) -\underline{u}\sin(\pi),$ where $\underline{u}$ is a pure unit quaternion. A simple replacement does not work in this case due to the noncommutativity of multiplication of quaternions. 
But we can use the multiplication from the righthand-side 
$$ R(q)f(\underline{x})=f(\underline{x})q,\quad f\in L^p(\mathbb{R}^3,\mathbb{H}). $$

\begin{defn}[Quaternionic fractional Riesz-Hilbert transform] For $f\in L^p(\mathbb{R}^3,\,\mathbb{H})$  the quaternionic fractional Riesz-Hilbert transform is defined as
\begin{eqnarray} \mathcal{H}^{\underline{u}\alpha} := \cos(\tfrac{\pi}{2}\alpha) + \sin(\tfrac{\pi}{2} \alpha) R(\underline{u})H \quad \mbox{and}\quad  H^{\underline{u}\alpha} := e^{-R(\underline{u})\tfrac{\pi}{2}\alpha}\mathcal{H}^{\underline{u}\alpha}. \label{QFRHT}
\end{eqnarray}
\end{defn}
There is a huge difference between formula (\ref{FHRT}) and (\ref{QFRHT}). The first formula is essentially a rotation in the plane spanned by $f$ and $Hf$, where as the second formula describes a rotation in $\mathbb{R}^4.$
\begin{thm}
The fractional Riesz-Hilbert transforms $H^{\underline{u}\alpha},\,\mathcal{H}^{\underline{u}\alpha}: L^p(\mathbb{R}^3,\,\mathbb{H}) \to L^p(\mathbb{R}^3,\,\mathbb{H}),\ 1<p<\infty, \,\alpha\in\mathbb{R},$ are linear, shift and scale invariant and fulfill the following properties
\begin{enumerate}
\item the inverse of $H^{\underline{u}\alpha}$ is $H^{-\underline{u}\alpha}$ and the inverse of $\mathcal{H}^{\underline{u}\alpha}$ is $\mathcal{H}^{-\underline{u}\alpha}, \alpha \in \mathbb{R},$ 
\item $H^{\underline{u}\alpha}$ is $2$-periodic in $\alpha,$ i.e. $H^{\underline{u}(\alpha +2)} = H^{\underline{u}\alpha},$ whereas $\mathcal{H}^{\underline{u}\alpha}$ is $4$-periodic in $\alpha,$ i.e. $\mathcal{H}^{\underline{u}(\alpha + 4)} = \mathcal{H}^{\underline{u}\alpha},\ \alpha\in\mathbb{R}.,$
\item If $f,\,g \in L^2(\mathbb{R}^3,\,\mathbb{R})$ such that $\langle f,\,g \rangle = 0$ then 
$$\langle H^{\underline{u}\alpha}f,\, H^{\underline{u}\alpha}g\rangle = \langle \mathcal{H}^{\underline{u}\alpha}f,\, \mathcal{H}^{\underline{u}\alpha}g\rangle = 0.$$
\end{enumerate}

\end{thm}
The proof of this theorem is similar to that for the fractional Riesz-Hilbert transform. We only prove the last property:
\begin{align*}
\langle & H^{\underline{u}\alpha}f,\, H^{\underline{u}\alpha}g\rangle =\int_{\mathbb{R}^3} \overline{(H^{\underline{u}\alpha}f)(\underline{x})}^{\mathbb{CH}}\,(H^{\underline{u}\alpha}f)(\underline{x})\,d\underline{x} \\ 
= & \int_{\mathbb{R}^3} e^{i\tfrac{\pi}{2}\alpha}\left(\cos(\tfrac{\pi}{2}\alpha)\overline{f(\underline{x})}^{\mathbb{CH}} - i\sin(\tfrac{\pi}{2}\alpha)\overline{(Hf)(\underline{x})}^{\mathbb{CH}}\right)e^{-i\tfrac{\pi}{2}\alpha}\left(\cos(\tfrac{\pi}{2}\alpha)g(\underline{x}) + i\sin(\tfrac{\pi}{2}\alpha)(Hg)(\underline{x})\right)\,d\underline{x} \\
= & \int_{\mathbb{R}^3} \overline{(\mathcal{H}^{\underline{u}\alpha}f)(\underline{x})}^{\mathcal{CH}}\,(\mathcal{H}^{\underline{u}\alpha}g)(\underline{x})\,d\underline{x} = \langle \mathcal{H}^{\underline{u}\alpha}f,\, \mathcal{H}^{\underline{u}\alpha}g\rangle \\
= & \cos^2(\tfrac{\pi}{2}\alpha)\langle f,\,g\rangle -i\sin(\tfrac{\pi}{2}\alpha)\cos(\tfrac{\pi}{2}\alpha)\left(\langle Hf,\,g\rangle -\langle f,\, Hg\rangle \right) + \sin^2(\tfrac{\pi}{2}\alpha)\langle Hf,\, Hg\rangle  = \langle f,\, g \rangle = 0.
\end{align*}
That last equation follows from Lemma 3.1 that implies $\langle Hf,\,g\rangle -\langle f,\, Hg\rangle = 0$ and $\langle Hf,\, Hg\rangle = \langle f,\,g\rangle $ (energy preservation).

\begin{thm} The quaternionic Riesz-Hilbert transform $R(\underline{u})H$ generates a norm continuous $C_0$-semigroup in $L^p(\mathbb{R}^3,\mathbb{H}).$ We have
$$ e^{\frac{\pi}{2}\alpha\,R(\underline{u})\, H} =\cos\left(\frac{\pi}{2}\alpha\right)I + \sin\left(\frac{\pi}{2}\alpha\right)R(\underline{u}) H= \mathcal{H}^{\underline{u}\alpha} .$$
The operator $R(\underline{u})(I-H)$ also generates a norm continuous $C_0$-semigroup in $L^p(\mathbb{R}^3,\mathbb{H}).$ We have
$$ e^{-\pi\alpha \frac{1}{2}R(\underline{u})(I-H)} 
     = R(e^{-\underline{u}\frac{\pi}{2}\alpha})\left(\cos\left(\frac{\pi}{2}\alpha\right) I + \sin\left(\frac{\pi}{2}\alpha\right) R(\underline{u})H \right)  = H^{\underline{u}\alpha}. $$
\end{thm}
Proof: Due to $(R(\underline{u})H)^2f = R(\underline{u})H(Hf\underline{u}) = (H^2f)\underline{u}^2 = -f $ we have $(R(\underline{u})H)^2 = -I$ and hence
\begin{align*}
e^{\tfrac{\pi}{2}\alpha R(\underline{u})H} & = \sum_{l=0}^{\infty} \frac{\left(\tfrac{\pi}{2}\alpha\right)^{2l}}{(2l)!} \left(R(\underline{u})H\right)^{2l} + \sum_{l=0}^{\infty} \frac{\left(\tfrac{\pi}{2}\alpha\right)^{2l+1}}{(2l+1)!} \left(R(\underline{u})H\right)^{2l+1} \\
  & = \sum_{l=0}^{\infty} \frac{\left(\tfrac{\pi}{2}\alpha\right)^{2l}}{(2l)!} \left(-1\right)^{l} + \sum_{l=0}^{\infty} \frac{\left(\tfrac{\pi}{2}\alpha\right)^{2l+1}}{(2l+1)!} (-1)^lR(\underline{u})H \\
  & = \cos(\tfrac{\pi}{2}\alpha) + \sin(\tfrac{\pi}{2} \alpha) R(\underline{u})H  =  \mathcal{H}^{\underline{u}\alpha} 
\end{align*}
Further, because the operators $R(\underline{u})$ and $R(\underline{u})H$ due to $R(\underline{u})[R(\underline{u})H]f = R(\underline{u})[(Hf)\underline{u}] =
(Hf)\underline{u}^2 = R(\underline{u})[HR(\underline{u})]f = R(\underline{u})[H(f\underline{u})] = R(\underline{u})[(Hf)\underline{u}] = (Hf)\underline{u}^2$ commute, we obtain
$$ e^{-\tfrac{\pi}{2}\alpha R(\underline{u})} e^{\tfrac{\pi}{2}\alpha R(\underline{u})H} = e^{-\pi\alpha \tfrac{1}{2}R(\underline{u})(I-H)}. $$
\begin{align*}
\lim_{\alpha \to 0} ||\mathcal{H}^{\underline{u}\alpha}f - f||_{L^p} & = \lim_{\alpha\to 0} || \cos\left(\tfrac{\pi}{2}\alpha\right)If + i \sin\left(\tfrac{\pi}{2}\alpha\right)(Hf)\underline{u} - f ||_{L^p} \\
& \leq \lim_{\alpha\to 0} \left(|\cos(\tfrac{\pi}{2}\alpha) - 1| \,||f||_{L^p} + |\sin(\tfrac{\pi}{2}\alpha)|\, ||Hf||_{L^p}|\underline{u}|\right) = 0 \\
\mbox{und   } 
\lim_{\alpha \to 0} ||H^{\underline{u}\alpha}f - f||_{L^p} & = \lim_{\alpha \to 0} ||e^{-R(\underline{u})\tfrac{\pi}{2}\alpha}\mathcal{H}^{\underline{u}\alpha}f - f||_{L^p} \\
& = \lim_{\alpha \to 0} ||(e^{-R(\underline{u})\tfrac{\pi}{2}\alpha}-1)\mathcal{H}^{\underline{u}\alpha}f + \mathcal{H}^{\underline{u}\alpha}f - f||_{L^p} \\
& = \lim_{\alpha \to 0} || (\cos(\tfrac{\pi}{2}\alpha)-1)\mathcal{H}^{\underline{u}\alpha}f - \sin(\tfrac{\pi}{2}\alpha)(\mathcal{H}^{\underline{u}\alpha}f)\underline{u} +
\mathcal{H}^{\underline{u}\alpha}f - f ||_{L^p} \\
& \hspace*{-3cm} = \lim_{\alpha \to 0} \left( |\cos(\tfrac{\pi}{2}\alpha)-1|\,||\mathcal{H}^{\underline{u}\alpha}f ||_{L^p} + |\sin(\tfrac{\pi}{2}\alpha)|\, ||\mathcal{H}^{\underline{u}\alpha}f||_{L^p}|\underline{u}| + ||\mathcal{H}^{\underline{u}\alpha}f - f ||_{L^p}\right) = 0.
\end{align*}

\begin{defn}[Quaternionic fractional monogenic signal] For $f\in L^p(\mathbb{R}^3,\mathbb{C})$ the fractional monogenic signal is defined as
$$ M^{\underline{u}\alpha}f := (I-R(e^{\underline{u}\frac{\pi}{2}\alpha})\mathcal{H}^{\underline{u}\alpha})f = (I+R(e^{\underline{u}\pi(1-\alpha)})H^{\alpha})f . $$
\end{defn}

\begin{thm}[Properties of the quaternionic fractional monogenic signal] 
Let $f\in L^p(\mathbb{R}^3,\mathbb{R}),$ then the quaternionic fractional monogenic signal  is linear, shift- and scale-invariant and possesses the following properties
\begin{enumerate}
\item $M^{\underline{u}\alpha}f \in H^+(\mathbb{R}^3),$ i.e. boundary values of a left-monogenic function in the upper half space.
\item $||M^{\underline{u}\alpha}f(\underline{x})||_{L^p} = |\sin(\frac{\pi}{2}\alpha)|\, ||Mf||_{L^p}.$ 
\item The quaternionic fractional monogenic signal is a rotated and modulated version of the monogenic signal.
\item The quaternionic fractional monogenic signal is not orthogonal to the signal and
\begin{align*}
Mf  & = (I+H)f = (1- \cot (\tfrac{\pi}{2}\alpha)C_r(\underline{u}))f - \csc (\tfrac{\pi}{2}\alpha) C_r(e^{\underline{u}\alpha\tfrac{\pi}{2}}\underline{u})H^{\underline{u}\alpha}f \\ 
      & = f - \cot (\tfrac{\pi}{2}\alpha)f\underline{u} - \csc (\tfrac{\pi}{2}\alpha) (H^{\underline{u}\alpha}f) e^{\underline{u}\alpha\tfrac{\pi}{2}}\underline{u} 
\end{align*}
\item The quaternionic fractional monogenic signal is $2$-periodic.
\end{enumerate}
\end{thm}
Proof: 
\begin{align*}
M^{\underline{u}\alpha}f  & = (I-R(e^{\underline{u}\frac{\pi}{2}\alpha})\mathcal{H}^{\underline{u}\alpha})f = (I-\left(\cos(\tfrac{\pi}{2}\alpha) + \sin(\tfrac{\pi}{2} \alpha) R(\underline{u})H \right)e^{\underline{u}\frac{\pi}{2}\alpha})f \\
 & = (I-\left(\cos(\tfrac{\pi}{2}\alpha) + \sin(\tfrac{\pi}{2} \alpha) R(\underline{u})H \right)\left(\cos(\tfrac{\pi}{2}) + \underline{u} \sin(\tfrac{\pi}{2}\alpha)\right))f \\
 & = (I-\cos^2(\tfrac{\pi}{2}\alpha) -\sin(\tfrac{\pi}{2}\alpha)\cos(\tfrac{\pi}{2}\alpha)H\underline{u} -\sin(\tfrac{\pi}{2}\alpha)\cos(\tfrac{\pi}{2}\alpha)\underline{u} -\sin^2(\tfrac{\pi}{2}\alpha)H\underline{u}^2)f \\
 & = \sin(\tfrac{\pi}{2}\alpha)(I+H)(\sin(\tfrac{\pi}{2}\alpha) - \cos(\tfrac{\pi}{2})\underline{u} )f  = ((I+H)f)\,e^{-\underline{u}\tfrac{\pi}{2}\alpha}(-\underline{u})\sin(\tfrac{\pi}{2}\alpha) \\
 & =  ((I+H)f)\,(\sin^2(\tfrac{\pi}{2}\alpha) - \sin(\tfrac{\pi}{2}\alpha)\cos(\tfrac{\pi}{2}\alpha)\underline{u})(-\underline{u})\\
 & =  ((I+H)f)\,\tfrac{1}{2}(-1+\cos(\pi\alpha) + \sin(\pi\alpha)\underline{u})\underline{u}
\end{align*}
and the quaternioinc fractional monogenic signal is 2-periodic in $\alpha$ and further
$$ ||M^{\underline{u}\alpha}f ||_{L^p} = |\sin(\tfrac{\pi}{2}\alpha)|\,||Mf||_{L^p} . $$
Because the fractional Riesz-Hilbert transform can be expressed by the Riesz-Hilbert transform we get
\begin{align*}
Mf & = (I+H)f  
= \left(f - \left((H^{\underline{u}\alpha}f)\,e^{\underline{u}\alpha\tfrac{\pi}{2}} - \cos(\alpha\tfrac{\pi}{2})f\right) \frac{\underline{u}}{\sin(\tfrac{\pi}{2}\alpha)}\right) \\
& = \left(I-\cot(\tfrac{\pi}{2}\alpha)R(\underline{u}) - \csc(\tfrac{\pi}{2})R(e^{\underline{u}\alpha\tfrac{\pi}{2}}\underline{u})H^{\underline{u}\alpha}\right)f \\
& = f - \cot (\tfrac{\pi}{2}\alpha)f\underline{u} - \csc (\tfrac{\pi}{2}\alpha) (H^{\underline{u}\alpha}f) e^{\underline{u}\alpha\tfrac{\pi}{2}}\underline{u}.
\end{align*}


\begin{thebibliography}{10}
\bibitem{Ball} R.S. Ball, ed. H. Gravelius,\textit{ Theoretische Mechanik starrer Systeme. Auf Grund der Methoden und Arbeiten mit einem Vorworte.} Berlin: Georg Reimer, 1889.
\bibitem{BBRH} S. Bernstein, J.-L. Bouchot, M. Reinhardt, B. Heise, \textit{Generalized Analytic Signals in Image Processing: Comparison, Theory and Applications,} in: E. Hitzer and S.J. Sangwine (eds), Quaternion and Clifford Fourier Transforms and Wavelets, Trends in Mathematics, Birkh\"auser, (2013), 221--246.
\bibitem{BE1} S. Bernstein, \textit{A Fractional Hilbert Transform for 2D Signals,} Adv. Appl. Clifford Algebras 24 (2014), 921--930.
\bibitem{BE2} S. Bernstein, \textit{The Fractional Monogenic Signal,} in S. Bernstein et al ed. Hypercomplex Analysis: New Perspectives and Applications, TIM, Birkhäuser Basel, (2014), 75--88.
\bibitem{Cu} A. Cusmariu, \textit{Fractional analytic signals.} Signal Processing \textbf{82} (2002), 267--272.
\bibitem{DNC} J.A. Davis, D.E. McNamara, D.M. Cottrell, \textit{Analysis of the fractional Hilbert transform.} Appl. Optics, \textbf{37} (1998), 6911--6913.
\bibitem{De} R. Delanghe, \textit{Clifford Analysis: History and Perspective.} Comp. Meth. Func. Theory, \textbf{1}(1) (2001), 107--153.
\bibitem{FS} M. Felsberg, G. Sommer, \textit{The monogenic signal.} IEEE Trans. Signal Proc., \textbf{49}(12) (2001), 3136--3144.
\bibitem{Ga} D. Gabor, \textit{Theory of communication.} J. of the Institution of Electrical Engineers - Part III:
Radio and Communication Engineering, \textbf{93}(26) (1946), 429--457.
\bibitem{Ha} S.L. Hahn, \textit{Multidimensional complex signals with single-orthant spectra.} Proc. IEEE, \textbf{80}(8) (1992), 1287--1300.
\bibitem{Kastler} A. Kastler, Rev. Opt. \textbf{29}, 308 (1950).
\bibitem{Larkin} K.G. Larkin, D.J. Bone, M.A. Oldfield, \textit{Natural demodulation of two-dimensional fringe patterns. I. General background of the spiral phase quadrature transform,} J. Optical Society of America A,  \textbf{18}(8), (2001), 1862--1870.
\bibitem{LMZ} A.W. Lohmann, D. Mendlovic, Z. Zalevsky, \textit{Fractional Hilbert transform.} Optics Letters, \textbf{21} (1996), 281--283.
\bibitem{Lu} P. Lounesto, \textit{Clifford Algebras and Spinors,} Cambridge Univ. Press, 1997.
\bibitem{Mebi} J.E. Mebius, \textit{A Matrix-based Proof of the Quaternion Representation Theorem for Four-Dimensional Rotations}, http://arxiv.org/abs/math/0501249v1.
\bibitem{McI} A. McIntosh, \textit{Fourier theory, singular integrals and harmonic functions on Lipschitz domains,} in: J. Ryan (ed.), Clifford Algebras in Analysis and Related Topics, CRC Press, (1996), 33--88.
\bibitem{Ryan} T. Qian, J. Ryan, \textit{Conformal Transformations and Hardy spaces arising in Clifford analysis,} J. Operator theory, \textbf{35}(1996), 349--372.
\bibitem{V} K. Th. Vahlen, Über Bewegungen und Complexe Zahlen, Math. Ann. 55, 1902,
 585-593.
\bibitem{VS} A. Venkitaraman, C.S. Seelamantula, \textit{Fractional Hilbert transform extensions and associated analytic signal construction.} Signal Processing, \textbf{94} (2014), 359--372.
 \bibitem{WW} J.L. Weiner, G.R. Wilkens, \textit{Quaternions and Rotations in $\mathbb{E}^4$.} Amer. Math. Monthly, \textbf{112}(1) (2005), 69--76.
\bibitem{Wolter} H. Wolter, Ann. Phys., \textbf{7}, 341 (1951).
\bibitem{Za} A.I. Zayed, \textit{Hilbert transform associated with the fractional Fourier transform.} IEEE Signal Processing Letters, \textbf{5} (1998), 206--208.
\end{thebibliography}
\end{document}